\documentclass[12pt]{article} 
\usepackage{amsfonts,latexsym}
\newlength\cellsize 
\savebox2{%
\begin{picture}(18,18)
\put(0,0){\line(1,0){18}}
\put(0,0){\line(0,1){18}}
\put(18,0){\line(0,1){18}}
\put(0,18){\line(1,0){18}}
\end{picture}}
\newcommand\cellify[1]{\def\thearg{#1}\def\nothing{}%
\ifx\thearg\nothing
\vrule width0pt height\cellsize depth0pt\else
\hbox to 0pt{\usebox2\hss}\fi%
\vbox to 18\unitlength{
\vss
\hbox to 18\unitlength{\hss$#1$\hss}
\vss}}
\newcommand\tableau[1]{\vtop{\let\\=\cr
\setlength\baselineskip{-16000pt}
\setlength\lineskiplimit{16000pt}
\setlength\lineskip{0pt}
\halign{&\cellify{##}\cr#1\crcr}}}
\savebox3{%
\begin{picture}(15,15)
\put(0,0){\line(1,0){15}}
\put(0,0){\line(0,1){15}}
\put(15,0){\line(0,1){15}}
\put(0,15){\line(1,0){15}}
\end{picture}}
\newcommand\expath[1]{%
\hbox to 0pt{\usebox3\hss}%
\vbox to 15\unitlength{
\vss
\hbox to 15\unitlength{\hss$#1$\hss}
\vss}}

\begin{document}
\newcommand{\qbc}[2]{ {\left [{#1 \atop #2}\right ]}}
\newcommand{\anbc}[2]{{\left\langle {#1 \atop #2} \right\rangle}}
\newcommand{\be}{\begin{enumerate}}
\newcommand{\ee}{\end{enumerate}}
\newcommand{\beq}{\begin{equation}}
\newcommand{\eeq}{\end{equation}}
\newcommand{\bea}{\begin{eqnarray}}
\newcommand{\eea}{\end{eqnarray}}
\newcommand{\beas}{\begin{eqnarray*}}
\newcommand{\eeas}{\end{eqnarray*}}
\newcommand{\fs}{\mathfrak{S}}
\newcommand{\cq}{{\cal Q}}
\newcommand{\im}{\mathrm{Imm}}
\newcommand{\bp}{\begin{problem}$\!\!\!$\textbf{.}\ }
\newcommand{\bt}{\begin{theorem}$\!\!\!$\textbf{.}\ }
\newcommand{\ep}{\end{problem}}
\newcommand{\et}{\end{theorem}}
\newcommand{\ds}{\displaystyle}
\newcommand{\sn}{\mathfrak{S}_n}
\newcommand{\wt}{\mathrm{wt}}
\newcommand{\co}{\mathrm{co}}
\newcommand{\cow}{\mathrm{co}(w^{-1})}
\newcommand{\bm}[1]{{\mbox{\boldmath $#1$}}}
\newcommand{\sbm}[1]{{\mbox{\scriptsize\boldmath $#1$}}}
\newcommand{\zz}{\mathbb{Z}}
\newcommand{\qq}{\mathbb{Q}}
\newcommand{\rr}{\mathbb{R}}
\newcommand{\pp}{\mathbb{P}}
\newcommand{\nn}{\mathbb{N}}
\newcommand{\cc}{\mathbb{C}}
\newcommand{\comp}{\mathrm{Comp}}
\newcommand{\st}{\,:\,}
\newcommand{\inv}{\mathrm{inv}}
\newcommand{\maj}{\mathrm{maj}}
\newcommand{\des}{\mathrm{des}}
\newcommand{\pr}{\mathrm{Prob}}

\newcounter{environment}[section]
\renewcommand{\theenvironment}{%
\arabic{section}.\arabic{environment}}

\newenvironment{definition}%
{\begin{rm}\refstepcounter{environment}{\bf\theenvironment\
Definition.~~}}%
{\end{rm}}

\newenvironment{example}%
{\begin{rm}\refstepcounter{environment}{\bf\theenvironment\
Example.~~}}%
{\end{rm}}

\newenvironment{proposition}%
{\begin{rm}\refstepcounter{environment}{\bf\theenvironment\
Proposition.~~}}%
{\end{rm}}

\newenvironment{propositiondefinition}%
{\begin{rm}\refstepcounter{environment}{\bf\theenvironment\
Proposition and Definition.~~}}%
{\end{rm}}

\newenvironment{theorem}%
{\begin{rm}\refstepcounter{environment}{\bf\theenvironment\
Theorem.~~}}%
{\end{rm}}

\newenvironment{corollary}%
{\begin{rm}\refstepcounter{environment}{\bf\theenvironment\
Corollary.~~}}%
{\end{rm}}

\newenvironment{lemma}%
{\begin{rm}\refstepcounter{environment}{\bf\theenvironment\
Lemma.~~}}%
{\end{rm}}

\newenvironment{conjecture}%
{\begin{rm}\refstepcounter{environment}{\bf\theenvironment\
Conjecture.~~}}%
{\end{rm}}

\begin{centering}
{\Large\bf {Generalized Riffle Shuffles}}\\
  {\Large\bf {and
  Quasisymmetric functions}}\\[.2in]
Richard P. Stanley\footnote{Partially supported by NSF grant
  DMS-9500714} \\
Department of Mathematics 2-375\\
Massachusetts Institute of Technology\\
Cambridge, MA 02139\\[.2in]
\textsc{Running title:} Generalized riffle shuffles and quasisymmetric
  functions\\[.5in] 
\end{centering}

\section{Introduction.}
\label{sec1}

\indent Let $x_i$ be a probability distribution on a totally ordered
set $I$, i.e., the probability of $i\in I$ is $x_i$. (Hence $x_i\geq
0$ and $\sum x_i=1$.) Fix $n\in\mathbb{P}=\{1,2,\dots\}$, and define a
random permutation $w\in\sn$ as follows. For each $1\leq j\leq n$,
choose independently an integer $i_j$ (from the distribution $x_i$).
Then \emph{standardize} the sequence $\bm{i}=i_1\cdots i_n$ in the
sense of \cite[p.~322]{ec2}, i.e., let $\alpha_1<\cdots<\alpha_k$ be
the elements of $I$ actually appearing in $\bm{i}$, and let $a_i$ be
the number of $\alpha_i$'s in $\bm{i}$. Replace the $\alpha_1$'s in
$\bm{i}$ by $1,2,\dots, a_1$ from left-to-right, then the $\alpha_2$'s
in $\bm{i}$ by $a_1+1, a_1+2,\dots,a_1+a_2$ from left-to-right, etc.
For instance, if $I=\mathbb{P}$ and $\bm{i}=311431$, then $w=412653$.
This defines a probability distribution on the symmetric group $\sn$,
which we call the \emph{QS-distribution} (because of the close
connection with quasisymmetric functions explained below). If we need
to be explicit about the parameters $x=(x_i)_{i\in I}$, then we will
refer to the \emph{QS$(x)$-distribution}.

The QS-distribution is not new, at least when $I$ is finite. It
appears for instance in \cite[pp.\ 
153--154]{d-f-p}\cite{fulman}\cite{i-t-w}\cite{kup}. Although these
authors anticipate some of our results, they don't make systematic use
of quasisymmetric and symmetric functions. Some additional work
using the viewpoint of our paper is due to Fulman \cite{fulman2}.

\begin{example}
As an example of a general result to be proved later
(Theorem~\ref{winsn}), let us compute $\pr(213)$, the probability that a
random permutation $w\in\mathfrak{S}_3$ (chosen from the
QS-distribution) is equal to the permutation
213. A sequence $\bm{i}=i_1i_2i_3$ will have the standardization 213
if and only if $i_2<i_1\leq i_3$. Hence
  $$ \pr(213) = \sum_{a<b\leq c}x_a x_b x_c. $$ 
\end{example}
\indent There is an alternative way to describe the QS-distribution.
Suppose we have a deck of $n$ cards. Cut the deck into
random packets of respective sizes $a_i>0,\ i\in I$, such that $\sum
a_i=n$ and the probability of $(a_i)_{i\in I}$ is 
  $$ \mathrm{Prob}\left( (a_i)_{i\in I}\right) = n!
      \prod_{i\in I} \frac{x_i^{a_i}}{a_i!}. $$
\noindent Then riffle shuffle these packets $P_i$ together 
into a pile $P$ in the manner described by Bayer and Diaconis
\cite{ba-di} (see also \cite{d-m-p}\cite{fulman});
namely, after placing $k$ cards on $P$, the probability that the next
card comes from the top of the current packet $P_j$ is
proportional to the current number of cards in 
$P_j$. This card is then placed at the bottom of $P$. The ordinary
dovetail or riffle shuffle \cite[p.~294]{ba-di} 
corresponds to the case $I=\{1,2\}$ and $x_1=x_2=1/2$. In this case
the original deck is cut according to the 
binomial distribution. More generally, if for fixed $q\in\mathbb{P}$
we cut the deck into $q$ packets (some possibly empty) according to
the $q$-multinomial distribution, then we obtain the $q$-shuffles of
Bayer and Diaconis \cite[p.~299]{ba-di} (where they use $a$ for our
$q$). The $q$-shuffle is identical to the QS-distribution for
$I=\{1,2\dots,q\}$ and $x_1=x_2=\cdots =x_q=1/q$. We will denote
this distribution by $U_q$. If we $q$-shuffle a deck and then
$r$-shuffle it, the distribution is the same as a single $qr$-shuffle
\cite[Lemma~1]{ba-di}. In other words, if $*$ denotes convolution of
probability distributions then $U_q*U_r=U_{qr}$. We extend this result
to the QS-distribution in Theorem~\ref{conv}.

In the next section we will establish the connection between the
QS-distribution and the theory of quasisymmetric functions. In
Section~\ref{sec3} we use known results from the theory of
quasisymmetric and symmetric functions to obtain results about the
QS-distribution, many of them direct generalizations of the work of
Bayer, Diaconis, and Fulman mentioned above. For instance, we show
that the probability that a random permutation $w\in\sn$ (always
assumed chosen from the QS-distribution) has $k$ inversions is equal
to the probability that $w$ has major index $k$. This is an analogue
of MacMahon's corresponding result for the \emph{uniform} distribution
on $\sn$ (see e.g.\ \cite[p.~23]{ec1}). A further result is that if $T$ is
a standard Young tableau of shape $\lambda$, then the probability that
$T$ is the insertion tableau of $w$ under the Robinson-Schensted-Knuth
algorithm is $s_\lambda(x)$, where $x=(x_i)_{i\in I}$ and $s_\lambda$
denotes a Schur function. 
% We use this result in Section~\ref{sec4} to
% compute the ``asymptotic shape'' as $n\rightarrow\infty$ of the
% partition $\lambda$ of $n$ that maximizes $s_\lambda(x)$ for fixed
% $x$. This result is much easier to obtain than the corresponding
% result for the uniform distribution on $\sn$ (see e.g.\ \cite[Exer.\
% 7.109(g)]{ec2} and the references cited there). 

\section{The QS-distribution and quasisymmetric functions}
\label{sec2}
In this section we show the connection between the QS-distribution and
quasisymmetric functions. In general our reference for quasisymmetric
functions and symmetric functions will be
\cite[Ch.~7]{ec2}. Quasisymmetric functions may be regarded as bearing
the same relation to compositions (ordered partitions of a nonnegative
integer) as symmetric functions bear to partitions. Namely, the
homogeneous symmetric functions of degree $n$ form a vector space of
dimension $p(n)$ (the number of partitions of $n$) and have many
natural bases indexed by the partitions of $n$. Similarly, the
homogeneous quasisymmetric functions of degree $n$ form a vector space
of dimension $2^{n-1}$ (the number of compositions of $n$) and have
many natural bases indexed by the compositions of $n$. 

A \emph{quasisymmetric function} $F(z)$ is a power series of bounded
degree, say with rational coefficients, in the variables
$z=(z_i)_{i\in I}$, with the following property: let 
$i_1<i_2<\cdots <i_n$ and $j_1<j_2<\cdots<j_n$, where $i_k,j_k\in
I$. If $a_1,\dots,a_k\in\pp$, then the coefficient of
$z_{i_1}^{a_1}\cdots z_{i_n}^{a_n}$ is equal to the coefficient of
$z_{j_1}^{a_1}\cdots z_{j_n}^{a_n}$ in $F(z)$. The set of all
quasisymmetric functions forms a graded $\qq$-algebra denoted $\cq$ .
The quasisymmetric functions that are homogeneous of degree $n$ form a
$\qq$-vector space, denoted $\cq_n$. If $|I|\geq n$ (including $|I|=
\infty$), then $\dim \cq_n=2^{n-1}$. 

Let $\alpha=(\alpha_1,\dots,\alpha_k)$ be a \emph{composition} of $n$,
i.e., $\alpha_i\in\pp$ and $\sum\alpha_i=n$. Let $\comp(n)$ denote the
set of all compositions of $n$, so $\#\comp(n)=2^{n-1}$ \cite[pp.\
14--15]{ec1}. Define 
  $$ S_\alpha =\{\alpha_1,\alpha_1+\alpha_2,\dots,\alpha_1+\cdots
       +\alpha_{k-1}\} \subseteq \{1,\dots,n-1\}. $$
The \emph{fundamental quasisymmetric function} $L_\alpha=L_\alpha(z)$
is defined by 
  $$ L_\alpha(z) = \sum_{{i_1\leq\cdots\leq i_n\atop 
      i_j<i_{j+1}\ \mbox{\scriptsize if}\ j\in S_\alpha}}
      z_{i_1}\cdots  z_{i_n}. $$
It is not hard to show that when $|I|\geq n$,
the set $\{L_\alpha\st\alpha\in\comp(n)\}$ is a $\qq$-basis for
$\cq_n$ \cite[Prop.~7.19.1]{ec2}. More generally, if $|I|=k$ then the
set $\{L_\alpha\st\alpha\in\comp(n),\ \mathrm{length}(\alpha)\leq k\}$
is a $\qq$-basis for $cq_n$, and $L_\alpha=0$ if length$(\alpha)>k$.

If $w=w_1w_2\cdots w_n\in \sn$ then let $D(w)$ denote the
\emph{descent set} of $w$, i.e.,
  $$ D(w) = \{ i\st w_i>w_{i+1}\}. $$
Write co$(w)$ for the unique composition of $n$ satisfying
$S_{\co(w)}= D(w)$. In other words, if
co$(w)=(\alpha_1,\dots,\alpha_k)$, then $w$ has descents exactly at 
$\alpha_1$, $\alpha_1+\alpha_2, \dots$, $\alpha_1+\cdots+\alpha_{k-1}$. 
For any composition $\alpha$ of $n$ it is easy to
see that the (possibly infinite) series $L_\alpha(x)$ (where as always
$x_i\geq 0$ and $\sum x_i=1$) is absolutely convergent and therefore
represents a nonnegative real number.

The main result on which all our other results depend is the
following.

\begin{theorem} \label{winsn}
\emph{Let $w\in \sn$. The probability  $\pr(w)$ that a permutation in
$\sn$ chosen from the QS-distribution is equal to $w$ is given by}
  $$ \pr(w) = L_{\co(w^{-1})}(x). $$
\end{theorem}
\indent
\textbf{Proof.} This result is equivalent to \cite[Lemma~3.2]{ge-r}
(see also \cite[Lemma~9.39]{reut}). Because 
of its importance here we present the (easy) proof. The integers $i$
for which $i+1$ appears somewhere to the left of $i$ in $w=w_1w_2\cdots
w_n$ are just the elements of $D(w^{-1})$. Let $i_1,\dots,i_n$ be
chosen independently from the distribution $x_i$. Let
$a_j=i_{w^{-1}(j)}$, i.e., if we write down $i_1\cdots i_n$ underneath
$w_1\cdots w_n$, then $a_j$ appears below $j$. In order for the
sequence $i_1\cdots i_n$ to standardize to $w$ it is necessary and
sufficient that $a_1\leq a_2\leq\cdots\leq a_n$ and that $a_j<a_{j+1}$
whenever $j\in D(w^{-1})$. Hence
  \beas \pr(w) & = & \sum_{{a_1\leq a_2\leq\cdots\leq a_n\atop
      a_i<a_{i+1}\ \mathrm{if}\ i\in D(w^{-1})}} x_{a_1}x_{a_2}\cdots
  x_{a_n}\\ & = & L_{\co(w^{-1})}(x).\ \ \Box \eeas

For some algebraic results closely related to Theorem~\ref{winsn},
see \cite{d-h-t}.

Two special cases of Theorem~\ref{winsn} are of particular interest.
The first is the $q$-shuffle distribution $U_q$. Let $\des(u)$ denote
the number of \emph{descents} of the permutation $u\in\sn$, i.e.,
  $$ \des(u) = \#\{i\st w(i)>w(i+1)\}. $$
It follows easily from Theorem~\ref{winsn} or by a direct argument
(see \cite[p.~364]{ec2}) that the probability
$\pr_{U_q}(w)$ that a random permutation in $\sn$ chosen from the
distribution $U_q$ is equal to $w$ is given by
  $$ \pr_{U_q}(w) = {q-\des(w^{-1})+n-1\choose n}q^{-n}. $$ 
This is just the probability of obtaining $w$ with a $q$-shuffle, as
defined in \cite[{\S}3]{ba-di}, confirming that
$U_q$ is just the $q$-shuffle distribution.

For fixed $n$, the uniform distribution $U$ on $\sn$ is given by
   \beq U=\lim_{q\rightarrow \infty} U_q. \label{eq:uuq} \eeq
Because of this formula many of our results below may be considered as
generalizations of the corresponding results for the uniform
distribution. These results for the uniform distribution are
well-known and classical theorems in the theory of permutation
enumeration. 

The second interesting special case of Theorem~\ref{winsn} is defined
by 
  $$ x_i=(1-t)t^{i-1}, $$ 
where $I=\pp$ and $0\leq t\leq 1$. Given $u=u_1u_2\cdots
u_n\in \sn$, define $\bar{u}=v_1v_2\cdots v_n$, where
$v_i=n+1-u_{n+1-i}$. (Equivalently, $\bar{u}=w_0uw_0$, where
$w_0(i)=n+1-i$.) Set $e(u)=\mathrm{maj}(\bar{u})$, where maj denotes
\emph{major index} (or \emph{greater index}) of $u$ \cite[p.\
216]{ec1}, i.e.,  
  $$ \maj(u) = \sum_{i\in D(u)}i. $$
It follows from \cite[Lemma~7.19.10]{ec2} that
  \beq \pr(w) = \frac{t^{e(w^{-1})}(1-t)^n}{(1-t)(1-t^2)\cdots(1-t^n)}. 
    \label{majsp} \eeq
\indent The QS-distribution defines a Markov chain (or random walk)
  ${\cal M}_n$ on $\sn$ by defining the 
probability $\pr(u,uw)$ of the transition from $u$ to $uw$ to be $\pr(w) =
L_{\co(w^{-1})}(x)$. The next result determines the eigenvalues of
this Markov chain. Equivalently, define 
  $$ \Gamma_n(x) = \sum_{w\in\sn} L_{\co(w^{-1})}(x)w\in\mathbb{R}\sn,
  $$
where $\mathbb{R}\sn$ denotes the group algebra of $\sn$ over
$\mathbb{R}$. Then the eigenvalues of ${\cal M}_n$ are just the
eigenvalues of $\Gamma_n(x)$ acting on $\mathbb{R}\sn$ by right
multiplication. 

\begin{theorem} \label{eig}
\emph{The eigenvalues of ${\cal M}_n$ are the power sum symmetric
  functions 
$p_\lambda(x)$ for $\lambda\vdash n$. The eigenvalue $p_\lambda(x)$
occurs with multiplicity $n!/z_\lambda$, the number of elements in
$\sn$ of cycle type $\lambda$.}\\ 
\end{theorem}

\vspace{-1em}
The above theorem appears implicitly in \cite[Thm.~4.4]{g-r} and more
explicitly in \cite[Note~5.20]{ncsf2}. See also \cite[{\S}7]{ncsf3} and
\cite[Exer.\ 7.96]{ec2}. It is also a special case of
\cite[Thm.~1.2]{b-h-r}, as we now explain. 
% For the special case of the $U_q$-distribution, a proof 
% appears in \cite[??]{ba-di}.

Bidigare, Hanlon, and Rockmore \cite{b-h-r} define a random walk on
the set ${\cal R}$ of regions of a (finite) hyperplane arrangement
${\cal A}$ in $\mathbb{R}^n$ as follows. (Precise definitions of the 
terms used below related to arrangements may be found in \cite{b-h-r}.
For further information see \cite{br-d}.)  Let ${\cal F}$ be the set of
\emph{faces} of ${\cal A}$, i.e., the nonempty faces of the closure of
the regions of ${\cal A}$. Let wt be a probability distribution on
${\cal F}$. Given $R\in {\cal R}$ and $F\in {\cal F}$, define $F\cdot
R$ to be that region $R'$ of ${\cal A}$ closest to $R$ that has $F$ as
a face.  (There is always a unique such region $R'$.) Now given a
region $R$, choose $F\in {\cal F}$ with probability wt$(F)$ and move
to the region $F\cdot R$.

Consider the special case when ${\cal A}$ is the braid arrangement
${\cal B}_n$, i.e., the set of hyperplanes $u_i=u_j$ for $1\leq
i<j\leq n$, where $u_1,\dots,u_n$ are the coordinates in
$\mathbb{R}^n$. The set of regions of ${\cal B}_n$ can be identified
in a natural way with $\sn$, viz., identify the region $u_{a_1} <
u_{a_2} < \cdots<u_{a_n}$ with the permutation $w$ given by $w(a_i)
=i$. The faces of ${\cal B}_n$ correspond to \emph{ordered set
partitions} $\pi$ of $1,2,\dots,n$, i.e., $\pi=(B_1,\dots,B_k)$ where
$B_i\neq \emptyset$, $B_i\cap B_j=\emptyset$ for $i\neq j$, and $\cup
B_i =\{1,2,\dots,n\}$; viz., $\pi$ corresponds to the face defined by
$u_r=u_s$ if $r,s\in B_i$ for some $i$, and $u_r<u_s$ if $r\in B_i$
and $s\in B_j$ with $i<j$. Define the \emph{type} of $\pi$ to be the
composition 
  $$\mathrm{type}(\pi)=(\#B_1, \dots,\#B_k) \in\comp(n). $$ 
Given $\alpha=(\alpha_1,\dots,\alpha_k) \in \comp(n)$, define the
\emph{monomial quasisymmetric function} $M_\alpha(z)$ by
  $$ M_\alpha(z) = \sum_{i_1<\cdots<i_k}z_{i_1}^{\alpha_1}\cdots
         z_{i_k}^{\alpha_k}. $$
(See \cite[(7.87)]{ec2}.)  Now if $F$ is a face of ${\cal B}_n$
corresponding to the ordered partition $\pi$, then set
   $$ \mathrm{wt}(F) = M_{\mathrm{type}(\pi)}(x). $$
It is easy to see that wt is a probability distribution on the set of
faces of ${\cal B}_n$; in fact, for indeterminates $(z_i)_{i\in I}$
we have 
  $$ \sum_\pi M_{\mathrm{type}(\pi)}(z) = 
   \left( \sum z_i\right)^n. $$
Thus wt defines a Bidigare-Hanlon-Rockmore random walk on the regions
of ${\cal B}_n$. Let $P(R,R')$ denote the probability of moving from a
region $R$ to a region $R'$.

\begin{theorem} \label{thm:bhr}
\emph{Let $R$ correspond to $u\in \sn$ and $R'$ to $uw$. Then} 
  $$ P(R,R')=L_{\co(w^{-1})}(x). $$
\emph{Hence the Bidigare-Hanlon-Rockmore random walk just defined on
  the regions of ${\cal B}_n$ is
isomorphic to the random walk on $\sn$ defined by the QS-distribution.}
\end{theorem}

\textbf{Proof.} By symmetry we may suppose that $u=12\cdots n$, the
identity permutation. The faces $F$ of $R'$ for which $F\cdot R'=R$
correspond to those ordered partitions $\pi$ obtained from
$w^{-1}=w_1\cdots w_n$ by drawing bars between certain pairs $w_i$ and
$w_{i+1}$, where there must be a bar when $i\in D(w^{-1})$. The
elements contained between the consecutive bars (including a bar
before $w_1$ and after $w_n$) form the blocks of $\pi$, read from
left-to-right. For instance, if $w=582679134$, then one way of drawing
bars is $|58|267|9|1|34|$, yielding the ordered partition $\pi =
(\{5,8\},\{2,6,7\},\{9\},\{1\},\{3,4\})$. It follows that
  $$ P(R,R') =\sum_\beta M_\beta(x), $$
where $\beta$ runs over all compositions refining $\co(w^{-1})$. From
the definition of $L_\alpha(x)$ we have
   $$ \sum_\beta M_\beta(x) = L_{\co(w^{-1})}(x), $$
and the proof follows. $\ \ \Box$

The next result determines the convolution of the QS$(x)$-distribution
with the QS$(y)$-distribution. We write $xy$ for the variables
$x_iy_j$ in the order $x_iy_j<x_ry_s$ if either $i<r$ or $i=r$ and
$j<s$. 

\begin{theorem} \label{conv}
\emph{Suppose that a permutation $u\in\sn$ is chosen from the
QS$(x)$-distribution, and a permutation $v\in\sn$ is chosen from the
QS$(y)$-distribution. Let $w\in \sn$. Then the probability that $uv=w$
is equal to $L_{\co(w^{-1})}(xy)$. In other words, if $\ast$ denotes
convolution then}
  $$ \mathrm{QS}(x)\ast \mathrm{QS}(y) =\mathrm{QS}(xy). $$
\emph{Equivalently, in the ring $\mathbb{R}\sn$ we have}
  $$ \Gamma_n(x)\Gamma_n(y) = \Gamma_n(xy). $$
\end{theorem}
\indent
Theorem~\ref{conv} is equivalent to
\cite[(62)]{ncsf2}\cite[Thm.~9.37]{reut}. For the $U_q$-distribution
it was also proved in \cite[Lemma~1]{ba-di}, and this proof can be
easily extended to prove Theorem~\ref{conv} in its full generality.

\section{Enumerative properties of the QS-distribution}
\label{sec3}

Suppose that $X$ is any subset of $\sn$. Let $w\in\sn$ be chosen
from the QS-distribution. It follows from Theorem~\ref{winsn} that the
probability that $w\in X$ is given by
  \beq \pr(w\in X) = \sum_{u\in X}L_{\co(u^{-1})}(x). \label{eq:pwinx}
  \eeq 
There are many known instances of the set $X\subseteq \sn$ for which
the right-hand side of (\ref{eq:pwinx}) can be explicitly
evaluated. Most of this section will be devoted to some examples of
this type.  

Our first result involves the symmetric function $l_n(z)$ defined by
  $$ l_n(z) = \frac 1n\sum_{d|n}\mu(d)p_d^{n/d}(z), $$
where $\mu$ denotes the usual number-theoretic M\"obius function. More
information on $l_n(z)$ may be found in \cite[Ch.~8]{reut} or
\cite[Exer.~7.88--7.89]{ec2}. 

\begin{theorem} \label{cyc}
  \emph{Let $w$ be a random permutation in $\sn$, chosen from the
  QS-distribution. The probability $\pr(\rho(w)=\lambda)$ that $w$ has
  cycle type $\lambda=\langle 1^{m_1} 2^{m_2}\cdots\rangle$ $\vdash n$
  (i.e., $m_i$ cycles of length $i$) is given by}
  $$ \pr(\rho(w)=\lambda)= \prod_{i\geq 1}h_{m_i}[l_i](x), $$
\emph{where brackets denote plethysm \cite[{\S}I.8]{macd}\cite[p.\
  447]{ec2}.} 
\end{theorem}

The proof of Theorem~\ref{cyc} follows from \cite[Thm.~8.23 and
Thm.~9.41(a)]{reut}. In the special case of the $U_q$-distribution,
the result appears in \cite[Thm.~A]{d-m-p}.

As above let maj$(w)$ denote the major index of $w\in\sn$, and let
inv$(w)$ denote the number of inversions of $w$  \cite[pp.\
20--21]{ec1}, i.e., 
  $$ \inv(w) = \#\{(i,j)\st 1\leq i<j\leq n,\ w(i)>w(j)\}. $$
Let $I_n(j)$ (respectively, $M_n(j)$) denote the
probability that $w\in\sn$ satisfies inv$(w)=j$ (respectively,
maj$(w)=j$) under the QS-distribution.

\begin{theorem} \label{im}
\emph{We have}
  \bea M_n(j) & = & I_n(j) \nonumber\\
   \sum_{n\geq 0}\sum_{j\geq 0} \frac{M_n(j)t^jz^n}
        {(1-t)(1-t^2)\cdots(1-t^n)} & = & \prod_{i\geq 1}
       \prod_{j\in I}
       \left( 1-t^{i-1}x_jz\right)^{-1}. \label{eq:cps}
%  \sum_{j\geq 0} M_n(j)t^j & = & \sum_{\lambda\vdash n} t^{b(\lambda)}
%     \frac{}{}
  \eea   
\end{theorem}
\indent
\textbf{Proof.} The standardization $w$ of a sequence
$\bm{i}=i_1i_2\cdots i_n$ satisfies $\maj(w)=\maj(\bm{i})$ and
$\inv(w) =\inv(\bm{i})$.  Hence by definition of the QS-distribution
we have
  \beas \sum_j M_n(j)t^j & = & \sum_{\sbm{i}=i_1\cdots i_n}
     t^{\maj(\sbm{i})} x_{i_1}\cdots x_{i_n}\\
    \sum_j I_n(j)t^j & = & \sum_{\sbm{i}=i_1\cdots i_n}
     t^{\inv(\sbm{i})} x_{i_1}\cdots x_{i_n}. \eeas
Thus if 
  \beas F_\lambda(t) & = & \sum_v t^{\maj(v)}\\
        G_\lambda(t) & = & \sum_v t^{\inv(v)}, \eeas
where $v$ ranges over all permutations of the multiset $\{
1^{\lambda_1}, 2^{\lambda_2},\dots\}$, then
  \beas \sum_j M_n(j)t^j & = & \sum_{\lambda\vdash n} F_\lambda(t)
       m_\lambda(x)\\
   \sum_j I_n(j)t^j & = & \sum_{\lambda\vdash n} G_\lambda(t)
       m_\lambda(x), \eeas
where $m_\lambda(x)$ denotes a monomial symmetric function.
But it is known \cite[(45), p.~97]{mem}\cite[Prop.~1.3.17]{ec2} that
  $$ F_\lambda(t) = G_\lambda(t) = \left[ {n\atop \lambda_1,
      \lambda_2,\dots}\right]_t, $$
a $q$-multinomial coefficient in the variable $t$. Hence $M_n(j) =
I_n(j)$. 

Now let $h_\lambda(x)$ denote a complete homogeneous symmetric
function \cite[{\S}7.5]{ec2}. It is easy to see \cite[Prop.\
7.8.3]{ec2} that
  $$ h_\lambda(1,t,t^2,\dots) = \prod_i \frac{1}{(1-t)(1-t^2)\cdots
      (1-t^{\lambda_i})}. $$
Hence
  $$ \sum_j M_n(j)t^j = (1-t)\cdots(1-t^n)\sum_{\lambda\vdash n}
      h_\lambda(1,t,t^2,\dots)m_\lambda(x). $$
Equation (\ref{eq:cps}) then follows immediately from the expansion
\cite[(7.10)]{ec2} of the Cauchy product
$\prod_{i,j}(1-x_iy_j)^{-1}$. Equation (\ref{eq:cps}) is also
equivalent to \cite[(78)]{ncsf1}.$\ \Box$ 

Other expansions of the Cauchy product lead to further formulas for
$\sum_j M_n(j)t^j$. In particular, from \cite[Prop.~7.20, Thm.\
7.12.1, and Cor.~7.21.3]{ec2} there follows (using notation from
\cite{ec2}) 
  \bea \sum_j M_n(j)t^j & = & \sum_{\lambda\vdash n} t^{b(\lambda)} 
   \frac{(1-t)\cdots(1-t^n)}{\prod_{u\in\lambda} \left( 1-t^{h(u)}
   \right)}s_\lambda(x) \nonumber\\ & = &
   \sum_{\lambda\vdash n} 
   \frac{(1-t)\cdots(1-t^n)}{\left( 1-t^{\lambda_1}\right)
    \cdots\left( 1-t^{\lambda_l}\right)}z_\lambda^{-1}p_\lambda(x).
   \label{eq:mnjp} \eea   
A famous result of MacMahon 
% \cite{??} 
asserts that 
  $$ \#\{ w\in\sn\st \maj(w)=j\} = \#\{ w\in\sn\st \inv(w)=j\}. $$
Since the uniform distribution $U$ on $\sn$ satisfies
$U=\lim_{q\rightarrow\infty}U_q$, it follows that Theorem~\ref{im} is
a generalization of MacMahon's result. In fact, the proof we have
given of Theorem~\ref{im} shows that it is equivalent to the
equidistribution of maj and inv on the set of permutations of any
finite \emph{multiset} (where the underlying set is linear
ordered). This result was also known, at least implicitly, to
MacMahon. For further references and some generalizations, see
\cite{b-w}. 

In principle one can use Theorem~\ref{im} to derive the moments of
$\maj(w)$ (or equivalently $\inv(w)$). Let $E$ denote expectation with
respect to the QS-distribution and $E_U$ expectation with respect
to the uniform distribution. The next result obtains the first two
moments of $\maj(w)$, stated for simplicity as the expectations of
$\maj(w)$ and $\maj(w)(\maj(w)-1)$ (instead of $\maj(w)^2$).

\begin{corollary} \label{moments}
\emph{Choose $w\in\sn$ from the QS-distribution. Then}
  \bea E(\maj(w)) & = & E_U(\maj(w)) -\frac 12{n\choose 2}
   p_2(x) \label{eq:emajw}\\  
   E(\maj(w)(\maj(w)-1)) & = & E_U(\maj(w)(\maj(w)-1))\\
   &  & \hspace{-6em}
   -3{n+1\choose 4}p_2(x) +\frac 43{n\choose 3}p_3(x) +
   \frac 32{n\choose 4}p_2(x)^2. \label{eq:e2majw} \eea
\end{corollary}
\textsc{Note.} It is easy to see (e.g.\ \cite[p.~16]{knuth3}) that
  \beas E_U(\maj(w)) & = & \frac 12{n\choose 2}\\
    E_U(\maj(w)(\maj(w)-1)) & = & \frac{n(n-1)(n-2)(9n+13)}{144}. \eeas
\textbf{Proof of Corollary~\ref{moments}.} Let
  $$ \Lambda_n(t) = \sum_jM_n(j)t^j. $$
Then
  $$ E(\maj(w)) = \Lambda'_n(1). $$
When we differentiate (\ref{eq:mnjp}) term-by-term and set $t=1$, the
only surviving terms will be from $\lambda=\langle 1^n\rangle$ and
$\lambda = \langle 21^{n-2}\rangle$. Hence
  $$ E(\maj(w)) = \frac{d}{dt} \left[ \frac{(1-t)\cdots (1-t^n)}
    {(1-t)^n}\frac{1}{n!}p_1(x)^n \right. $$
  \beq \qquad \left. +\frac{(1-t)\cdots (1-t^n)}{(1-t^2)(1-t)^{n-2}}
    \frac{1}{2\cdot(n-2)!}p_2(x)^n\right]_{t=1}. \label{eq:ddt} \eeq
Now $p_1(x)=\sum x_i=1$ and
  $$ \frac{d}{dt}\left[ \frac{(1-t)\cdots (1-t^n)}{(1-t)^n}
    \frac{1}{n!}\right] = E_U(\maj(w)). $$
It is routine to compute the second term on the right-hand side of
(\ref{eq:ddt}), obtaining equation (\ref{eq:emajw}).

A similar computation using $E(\maj(w)(\maj(w)-1)) = \Lambda''_n(1)$
yields (\ref{eq:e2majw}); we omit the details. $\ \Box$
 
It is clear from the above proof that for general $k\geq 1$,
$E(\mathrm{maj}(w)^k)$ is a linear combination, whose coefficients
are polynomials in $n$, of power sum symmetric functions
$p_\lambda(x)$ where $\lambda\vdash n$ and $\ell(\lambda)\geq
n-k$. Here $\ell(\lambda)$ denotes the length (number of parts) of
$\lambda$. 

We next consider the relationship between the QS-distribution and the
Robinson-Schensted-Knuth (RSK) algorithm \cite[{\S}7.11]{ec2}. Recall
that this algorithm associates a pair $(T,T')$ of standard Young
tableaux (SYT) of the same shape $\lambda\vdash n$ with a permutation
$w\in\sn$. We call $T$ the \emph{insertion tableau} of $w$, denoted
ins$(w)$. The shape $\lambda$ of $T$ is also called the \emph{shape}
of $w$, denoted sh$(w)=\lambda$.

\begin{theorem} \label{rsk}
  \emph{Choose $w\in\sn$ from the QS-distribution. Let $T$ be an SYT
  of shape $\lambda\vdash n$. Then the probability that} ins$(w)=T$
  \emph{is given by}
  $$ \pr(\mathrm{ins}(w)=T) = s_\lambda(x), $$ 
\emph{where $s_\lambda$ denotes a Schur function.}
\end{theorem}

\textbf{Proof.} Let
  $$ X_T =\{w\in\sn\st\mathrm{ins}(w)=T\}. $$ 
It follows from \cite[Thm.~7.19.2 and Lemma~7.23.1]{ec2} that
   $$ \sum_{w\in X_T} L_{\co(w^{-1})}(x) = s_\lambda(x). $$ 
The proof now follows from equation (\ref{eq:pwinx}). $\ \Box$

\begin{corollary} \label{cor:shla}
 \emph{Choose $w\in\sn$ from the QS-distribution, and let
  $\lambda\vdash n$. Then}
  $$ \pr(\mathrm{sh}(w) = \lambda)=f^\lambda s_\lambda(x), $$ 
 \emph{where
  $f^\lambda$ denotes the number of SYT of shape $\lambda$ (given
  explicitly by the Frame-Robinson-Thrall hook-length formula
  \cite[Cor.~7.21.6]{ec2}).}
\end{corollary}

\textbf{Proof.} There are $f^\lambda$ SYT's $T$ of shape $\lambda$,
and by Theorem~\ref{rsk} they all have probability $s_\lambda(x)$ of
being the insertion tableau of $w$. $\ \Box$

\textsc{Note.} The probability distribution Prob$(\lambda)=f^\lambda
s_\lambda(x)$ on the set of all partitions $\lambda$ of all
nonnegative integers is a specialization of the $z$-measure of Borodin
and Olshanski, as surveyed in \cite{b-o} (see also \cite{ok3}). 

It is easy to give an expression for the probability
that a permutation $w\in\sn$ chosen from the QS-distribution has a
fixed descent set $S\subseteq\{1,2,\dots,n-1\}$. In general, if
$\lambda$ and $\mu$ are partitions with $\mu\subseteq \lambda$ (i.e.,
$\mu_i\leq \lambda_i$ for all $i$), then $\lambda/\mu$ denotes the
\emph{skew shape} obtained by removing the diagram of $\mu$ from that
of $\lambda$ \cite[p.~309]{ec2}. One can then define (see
\cite[Def.~7.10.1]{ec2}) the \emph{skew Schur function}
$s_{\lambda/\mu}(x)$.  A \emph{border strip} (or \emph{rim hook} or
\emph{ribbon}) is a connected skew shape with no $2\times 2$ square
\cite[p.~245]{ec2}. Let $\alpha=(\alpha_1,\dots,\alpha_k)
\in\comp(n)$. Then there is a unique border strip $B_\alpha$ (up to
translation) with $\alpha_i$ squares in row $i$. For instance, the
border strip $B_{(3,2,1,2,4)}$ looks like
  $$ \tableau{ & & & & & \mbox{} & \mbox{} & \mbox{}\\
    & & & & \mbox{} & \mbox{}\\ & & & & \mbox{}\\ & & & 
    \mbox{} & \mbox{}\\ \mbox{} & \mbox{} & \mbox{} & \mbox{}}\ . $$
Some special properties of border strip Schur functions
$s_{B_\alpha}$ are discussed in \cite[{\S}7.23]{ec2}. 

\begin{theorem}
\emph{Let $w$ be a random permutation in $\sn$, chosen from the
  QS-distribution. Let $\alpha=(\alpha_1,\dots,\alpha_k) \in\comp(n)$,
  and recall that $S_\alpha=\{\alpha_1,\alpha_1+\alpha_2,
  \dots,\alpha_1+\cdots+  
  \alpha_{k-1}\}\subseteq \{1,2\dots,n-1\}$. Then the  
  probability that $w$ has descent set $S_\alpha$ is given by}
  $$ \pr(D(w)=S_\alpha) = s_{B_\alpha}(x). $$
\end{theorem}
\indent
\textbf{Proof.} Immediate from Theorem~\ref{winsn} and the fact
\cite[Cor.~7.23.4]{ec2} that
  $$ s_{B_\alpha} = \sum_{{w\in \sn\atop \alpha=\co(w)}}
  L_{\co(w^{-1})}. \qquad \Box $$

\section{Longest increasing subsequences.} \label{sec4}

It is a fundamental result of Schensted \cite{schen} (see also
\cite[Cor.~7.23.11]{ec2}) that if
sh$(w)=\lambda=(\lambda_1,\lambda_2,\dots)$, then $\lambda_1$ is the
length is$(w)$ of the longest increasing subsequence of $w$. (There is
a generalization due to Greene \cite{greene}\cite[Thm.~7.23.13 or
A1.1.1]{ec2} 
that gives a similar interpretation of any $\lambda_i$.) Thus the
expected value $E_U(n)$ of is$(w)$ for $w\in\sn$ under the uniform
distribution is given by 
  $$ E_U(n)=E_U(\mathrm{is}(w)) = \frac{1}{n!}\sum_{\lambda\vdash n}
       \lambda_1\left( f^\lambda\right)^2. $$
\noindent (See \cite[Exer.~7.109(b)]{ec2}.) It was shown by Vershik
and Kerov \cite{v-k} that $E_U(n)\sim 2\sqrt{n}$. Later Baik \emph{et
al.}\ \cite{b-d-j} obtained a much stronger result, viz., the exact
distribution (suitably normalized) of $\lambda_1$ (or equivalently
is$(w)$) in the limit $n\rightarrow \infty$. (This result was extended
to any $\lambda_i$ in \cite{b-o-o}\cite{ok}\cite{joh}.) We can ask
whether similar results hold for the QS-distribution. Many results
along these lines appear in \cite{biane}\cite[Ex.\
2]{b-ok}\cite{b-o-o}\cite{ok2}\cite{ok3}. Here we make some elementary
comments pointing in a different direction from these papers. 

It is clear from Corollary~\ref{cor:shla} that
  \beq E(\mathrm{is}(w)) = \sum_{\lambda\vdash n} \lambda_1
     f^\lambda s_\lambda(x), \label{eq:isw} \eeq
where as usual $E$ denotes expectation with respect to the
QS-distribution. This formula can be made more explicit in the special
case of the distribution $U_q$. Identify $\lambda\vdash n$ with its
Young diagram, and let $c(u)$ denote the content of the square $u\in
\lambda$ as defined in \cite[p.~373]{ec2}.

\begin{theorem} \label{thm:uqisw}
\emph{We have}
  \beq E_{U_q}(\mathrm{is}(w)) = \frac{1}{n!}\sum_{\lambda\vdash n}
        \lambda_1\left( f^\lambda\right)^2\prod_{u\in\lambda}
          \left( 1+q^{-1}c(u)\right). \label{eq:euq} \eeq
\end{theorem}
\indent
\textbf{Proof.} Let $s_\lambda(1^q)$ denote $s_\lambda(x)$ where $x_1=
\cdots=x_q=1$ and $x_i=0$ for $i>q$. It is well-known (e.g.\
\cite[Cor.~7.21.4 and 7.21.6]{ec2}) that
  $$ s_\lambda(1^q) = \frac{f^\lambda}{n!}\prod_{u\in\lambda}
        (q+c(u)). $$
Since $s_\lambda(x)$ is homogeneous of degree $n$, the proof follows
from equation (\ref{eq:isw}). $\ \Box$ 

Note that $E_{U_q}(\mathrm{is}(w))$ is a polynomial in $q^{-1}$, which
can be written
  \beq E_{U_q}(\mathrm{is}(w)) = E_U(\mathrm{is}(w)) + \frac{1}{n!}
      \sum_{\lambda\vdash n}\lambda_1\left( f^\lambda\right)^2
       \left( \sum_{u\in\lambda}c(u)\right)\frac 1q +\cdots. 
    \label{eq:euq1} \eeq
Let us mention that
  $$ \sum_{u\in\lambda}c(u) = \sum i(\lambda'_i-\lambda_i) =
         \sum {\lambda_i\choose 2}-\sum{\lambda'_i\choose 2}, $$
where $\lambda'=(\lambda'_1,\lambda'_2,\dots)$ denotes the conjugate
partition to $\lambda$. If we regard $n$ as fixed, then equation
(\ref{eq:euq1}) shows the rate at which $E_{U_q}(\mathrm{is}(w))$
converges to $E_U(\mathrm{is}(w))$ as $q\rightarrow \infty$. It is
therefore of interest to investigate the asymptotic behavior of the
coefficient  
  $$ F_1(n)=\frac{1}{n!}
  \sum_{\lambda\vdash n}\lambda_1\left( f^\lambda\right)^2
   \left( \sum_{u\in\lambda}c(u)\right) $$ 
as $n\rightarrow\infty$ (as well as the coefficient of $q^{-i}$ for
$i>1$). Numerical evidence suggests that $F_1(n)/n$ is a slowly
increasing function of $n$. The largest value 
of $n$ for which we have done the computation is $n=47$, giving
$F_1(47)/47\approx 0.6991$. Eric Rains suggests that
there is some reason to suspect that $F(n)/n$ grows like $n^{1/6}$. 

\textsc{Acknowledgment.} I am grateful to Jinho Baik, Fran\c{c}ois
Bergeron, Ken Brown, Persi Diaconis, Alain Lascoux, Andrei Okounkov,
and Craig Tracy for providing helpful suggestions and additional
references.


\begin{thebibliography}{99}

\bibitem{b-d-j} J. Baik, P. Deift, and K. Johansson, On the
  distribution of the length of the longest increasing subsequence of
  random permutations, J.\ Amer.\ Math.\ Soc., to appear.

\bibitem{ba-di} D. Bayer and P. Diaconis, Trailing the dovetail shuffle
  to its lair, Ann.\ Applied Probability \textbf{2} (1992),
  294--313. 

\bibitem{b-h-r} P. Bidigare, P. Hanlon, and D. Rockmore, A
  combinatorial description of the spectrum of the Tsetlin library and
  its generalization to hyperplane arrangements, Duke Math.\
  J.\ \textbf{99} (1999), 135--174.

\bibitem{biane} Ph.\ Biane, seminar talk at U.C.\ Berkeley, March
  2000. 

\bibitem{b-w} A. Bj\"orner and M. Wachs, $q$-hook length formulas for
  forests, J.\ Combinatorial Theory (A) \textbf{52} (1989),
  165--187.

\bibitem{b-ok} A. Borodin and A. Okounkov, A Fredholm determinant
  formula for Toeplitz determinants, Integral Equations and
  Operator Theory \textbf{37} (2000), 386--396, math.CA/9907165.

\bibitem{b-o} A. Borodin and G. Olshanski, $Z$-measures on partitions,
  Robinson-Schensted-Knuth correspondence, and $\beta=2$ random matrix
  ensembles, preprint, math.CO/9905189.

\bibitem{b-o-o} A. Borodin, A. Okounkov, and G. Olshanski, On
  asymptotics of Plancherel measures for symmetric groups, J.\ 
    Amer.\ Math.\ Soc.\ \textbf{13} (2000), 481--515,
  math.CO/9905032. 

\bibitem{br-d} K. Brown and P. Diaconis, Random walks and hyperplane
  arrangements, Ann.\ Probability \textbf{26} (1998),
  1813--1854. 

\bibitem{d-f-p} P. Diaconis, J. Fill, and J. Pitman, Analysis of top
  to random shuffles, Combinatorics, Probability, and
  Computing \textbf{1} (1992), 135--155. 

\bibitem{d-m-p} P. Diaconis, M. McGrath, and J. Pitman, Riffle
  shuffles, cycles, and descents, Combinatorica \textbf{15}
  (1995), 11--29.

\bibitem{d-h-t} G. Duchamp, F. Hivert, and J.-Y.\ Thibon, Une
  g\'en\'eralisation des fonctions quasi-sym\'etriques et des
  fonctions sym\'etriques non commutatives, C.\ R.\ Acad.\ Sci.\
  Paris S\'er.\ I Math.\  \textbf{328} (1999), 1113--1116. 

\bibitem{ncsf3} G. Duchamp, A. Klyachko, D. Krob, and J.-Y.\ Thibon,
Noncommutative symmetric functions III: Deformations of 
Cauchy and convolution algebras, Discrete Math.\ Theor.\
Comput.\ Sci.\ \textbf{1} (1997), 159--216. 

\bibitem{fulman} J. Fulman, The combinatorics of biased riffle
  shuffles, Combinatorica \textbf{18} (1998), 173--184.

\bibitem{fulman2} J. Fulman, Applications of symmetric functions to
  cycle and subsequence structure after shuffles, preprint,
  math.CO/0102176. 

\bibitem{g-r} A. Garsia and C. Reutenauer, A decomposition of
  Solomon's descent algebra, Advances in Math.\ \textbf{77}
  (1989), 189--262.

\bibitem{ncsf1} I. M. Gelfand, D. Krob, A. Lascoux, B. Leclerc,
  V. S. Retakh, and J.-Y.\ Thibon, Noncommutative symmetric
functions, Advances in Math.\ \textbf{112} (1995), 218--348.

\bibitem{ge-r} I. Gessel and C. Reutenauer, Counting permutations with
  given cycle structure and descent set, J.~Combinatorial Theory
  (A) \textbf{64} (1993), 189--215.

\bibitem{greene} C. Greene, An extension of Schensted's theorem,
  Advances in Math.\ \textbf{14} (1974), 254--265.

\bibitem{i-t-w} A. R. Its, C. A. Tracy, and H. Widom, Random words,
  Toeplitz determinants and integrable systems.\ I,
  preprint, math.CO/9909169.

\bibitem{joh} K. Johansson, Discrete orthogonal polynomial ensembles
  and the Plan\-cherel measure, preprint, math.CO/9906120.

\bibitem{knuth3} D. E. Knuth, The Art of Computer Programming,
  vol.\ 3, Sorting and Searching, second ed., Addison-Wesley,
  Reading, Massachusetts, 1998.

\bibitem{ncsf2} D. Krob, B. Leclerc, and J.-Y.\ Thibon, Noncommutative
  symmetric functions II: Tranformations of alphabets,
  Internat.\  J.\ Algebra Comput.\ \textbf{7} (1997), 181--264.

\bibitem{kup} G. Kuperberg, Random words, quantum statistics, central
  limits, random matrices, preprint, math.PR/9909104. 

\bibitem{lo-sh} B. F. Logan and L. A. Shepp, A variational problem for
  random Young tableaux, Advances in Math.\ \textbf{26} (1977),
  206--222. 

% \bibitem{l-t} M. Ledoux and M. Talagrand, Probability in Banach
% Spaces, Springer-Verlag, Berlin, 1991.

\bibitem{macd} I. G. Macdonald, Symmetric Functions and Hall
Polynomials, second ed., Oxford University Press, Oxford, 1995.

\bibitem{ok2} A. Okounkov, Infinite wedge and random partitions,
  preprint, math.RT/9907127.

\bibitem{ok} A. Okounkov, Random matrices and random permutations,
Internat.\ Math.\ Res.\ Notices \textbf{2000}, 1043--1095,
math.CO/9903176. 

\bibitem{ok3} A. Okounkov, SL$(2)$ and $z$-measures, preprint,
  math.RT/0002135. 

\bibitem{reut} C. Reutenauer, Free Lie Algebras, Oxford
  University Press, Oxford, 1993.

\bibitem{schen} C. E. Schensted, Longest increasing and decreasing
  subsequences, Canad.\ Math.\ J.\ \textbf{13} (1961),
  179--191. 

\bibitem{mem} R. Stanley, Ordered structures and partitions,
Mem.\ Amer.\ Math.\ Soc., no. 119 (1972).

\bibitem{ec1} R. Stanley, Enumerative Combinatorics, vol.\ 1,
  second printing, Cambridge University Press, New York/Cambridge,
  1996.

\bibitem{ec2} R. Stanley, Enumerative Combinatorics, vol.\ 2, 
Cambridge University Press, New York/Cambridge, 1999.

\bibitem{v-k} A. M. Vershik and S. V. Kerov, Asymptotic behavior of
    the Plancherel measure of the symmetric group and the limit form
    of Young tableaux, Soviet Math.\
    Dokl.\ \textbf{18} (1977), 527--531; translated from Dokl.\
    Akad.\ Nauk SSSR \textbf{233} (1977), 1024--1027. 

\end{thebibliography}
\end{document}